\documentclass[11pt]{article}

\usepackage{amsmath}
\usepackage{amssymb}
\usepackage{epsfig}

\def\per{\mathop{\rm per}\nolimits}
\def\det{\mathop{\rm det}\nolimits}

\def\DS{\mathop{\rm DS}\nolimits}
\def\per{\mathop{\rm per}\nolimits}
\def\tr{\mathop{\rm tr}\nolimits}
\def\dsd{\mathop{\rm DSD}\nolimits}
\def\Dev{\mathop{\rm Dev}\nolimits}
\def\one{{\bf 1}}
\def\dscf#1#2{#1^{\downarrow#2}}

\def\sharpsum{\mathop{\mathord{\sum}^\sharp}}
\def\ave{\mathop{\rm ave}\nolimits}
\def\avesharp{\mathop{\mathord{\ave}^\sharp}}

\def\norm#1{\Vert#1\Vert}
\def\normsq#1{\norm{#1}^2}

\def\Partition{{\cal P}}
\mathchardef\Real="323C

\def\vhalf{\textstyle{\frac12}}

\font\fiverm=cmr5
\def\lambdamax{\lambda_{\hbox{\fiverm max}}}
\def\gap{\mathop{\rm Gap}\nolimits}
\title{An asymptotic approximation for the permanent of a doubly stochastic matrix}
\author{Peter McCullagh \\ University of Chicago}

\begin{document}

\maketitle
\begin{abstract}
\noindent
A determinantal approximation is obtained for the permanent of a doubly
stochastic matrix.
For moderate-deviation matrix sequences, the asymptotic relative error
is of order $O(n^{-1})$.
\end{abstract}
\textsf{\textbf{keywords}: Doubly stochastic Dirichlet distribution;
Maximum-likelihood projection;
Sinkhorn projection}

\section{Introduction}
A non-negative matrix of order~$n$ with unit row and column sums is called
doubly stochastic.
To each strictly positive square matrix $Y$ there corresponds a unique re-scaled matrix
$A$ such that $Y_{ij} = n A_{ij} \exp(\beta_i + \gamma_j)$ where
$A$ is strictly positive with unit row and column sums.
The doubly stochastic projection $Y \mapsto A$ can be computed
by iterative proportional scaling
(Deming and Stephan, 1940),
although, in this context, it is usually called the Sinkhorn algorihm
(Sinkhorn, 1964;
Linial, Samordnitsky and Wigderson, 1998).
A non-negative matrix $Y$ containing some zero components is said to be scalable
if there exist real numbers $\beta_i, \gamma_j$
such that $Y_{ij} = n A_{ij} \exp(\alpha_i+\gamma_j)$
where $A$ is doubly stochastic.
In that case, the pattern of zeros in $Y$ is the same as the pattern of zeros in~$A$.
A non-negative matrix is scalable if and only if $\per(Y) > 0$.

The permanent of a square matrix $A$
\[
\per(A) = \sum_{\sigma} \prod_{i=1}^n A_{i, \sigma(i)}
\]
is the sum over permutations $\sigma\colon[n]\to[n]$ of products,
$A_{1,\sigma(1)}\cdots A_{n,\sigma(n)}$,
one component taken from each row and each column.
If $A_{ij} = \alpha_i \beta_j B_{ij}$, then the permanents are related by
\[
\per(A) = \per(B)\, \prod_{i=1}^n (\alpha_i \beta_j) .
\]
Thus, the problem of approximating the permanent of a scalable non-negative matrix
is reduced to the approximation of the corresponding doubly stochastic matrix,
which is the subject of this paper.

It is known that the permanent of a generic matrix is not computable
in polynomial time as a function of~$n$ (Valiant, 1979).
Since exact computation is not feasible for large matrices, most authors have sought
to approximate the permanent using Monte-Carlo methods.
For example,
Jerrum and Sinclair (1989) and Jerrum, Sinclair and Vigoda (2004)
use a Markov chain technique,
while Kou and McCullagh (2009) use an importance-sampling scheme
for $\alpha$-permanent approximation.
This paper develops an entirely different sort of computable approximation,
deterministic but asymptotic.

The determinantal formula (\ref{det_approx}) is an asymptotic approximation for the
permanent of a doubly stochastic matrix~$B$,
the conditions for which are ordinarily satisfied if the original matrix~$A$
has independent and identically distributed components.
The approximation is not universally valid for all matrix sequences,
but it is valid under moderate-deviation conditions, which can be checked.
For example, $A$~may be symmetric.

\section{Sinkhorn projection}
The scaling algorithm has an interpretation connected with maximum
likelihood for generalized linear models as follows.
Let the components of $Y$ be independent exponential variables
with means $\mu_{ij}$ such that $\log\mu_{ij}= \beta_i + \gamma_j$ lies in the
additive subspace {\sl row+col}.
The maximum-likelihood projection $Y\mapsto \hat\mu$ is such that
\[
\sum_i (Y_{ij}/\hat\mu_{ij} - 1) = \sum_j (Y_{ij}/\hat\mu_{ij} - 1) = 0
\]
(McCullagh and Nelder, 1989).
The residual matrix $Y/\hat\mu$ is strictly positive
with row and column sums equal to~$n$,
from which it follows that $A = n^{-1} Y/\hat\mu$ is doubly stochastic.
The residual deviance 
\[
%2\sum_{ij} Y_{ij}/\hat\mu_{ij}-1 - \log(Y_{ij}/\hat \mu_{ij} =
\Dev(A) = -2\sum_{ij} \log(n A_{ij})
\]
is a measure of total variability for strictly positive
doubly stochastic matrices, taking the value zero only for
the uniform matrix $A_{ij} = 1/n$.
The exponential assumption can be replaced by any gamma distribution
provided that the gamma index $\nu$ is held constant.
The re-scaled matrix $A$ thus generated is said to have
the doubly stochastic Dirichlet distribution $\dsd_n(\nu)$ with parameter~$\nu$.
The smaller the value of $\nu$, the more extreme the components of~$A$,
and the greater the deviance.
The distribution $\dsd_n(1)$ is not uniform with respect to Lebesgue
measure in the sense of Chatterjee, Diaconis and Sly (2010),
but it is presumably close for large~$n$.
The maximum-likelihood estimate $\hat\nu(A)$ is such that
\[
2\log(\hat\nu) - 2\psi(\hat\nu) = \Dev(A)/n^2,
\]
where $\psi$~is the derivative of the log gamma function.
It is apparent from simulations that if $A\sim\dsd_n(\nu)$,
\[
%2\hat\nu(A) \times (\log\per(A) - \log(n!)) = 1 + O_p(n^{-1})
2\hat\nu(A) \times \log(\per(nA)/n!) = 1 + O_p(n^{-1}),
\]
with variability decreasing in~$n$.
This limiting product appears to be invariably less than one for
non-Dirichlet matrices.

The scaling $Y \mapsto A$ can be accomplished either by iterative
weighted least squares, as is usually done for generalized linear models,
or by the Sinkhorn iterative proportional scaling algorithm.
For present purposes, the latter algorithm is preferred because
it is efficient and simple to implement.
Linial, Samorodnitsky and Wigderson (1998) provide a
a modification that guarantees convergence in polynomial time.
For moderate-deviation matrices, the time taken to compute the
Sinkhorn projection is usually negligible.

For a doubly stochastic matrix, it is known that 
$$
0 \le \log(\per(nA)/n!) \le n
$$
the lower bound being attained at the matrix $A=J$ whose components
are all equal,~$J_{ij} = 1/n$.
The upper bound of $\log(n^n/n!) \simeq n$ is attained at each of the permutation matrices,
which are the extreme points in the set $\DS_n$ of
doubly-stochastic matrices of order~$n$.
These bounds suggest that it may be easier to develop an approximation
for the permanent of a doubly stochastic matrix
than the permanent of an arbitrary positive matrix.

\section{Moderate-deviation sequences}
A sequence of square matrices $X_1, X_2,\ldots$ in which $X_n$ is of order~$n$
is called \emph{weakly bounded} if the absolute $p$th moment
\[
\mu_p = \limsup_{n\to\infty} \frac 1 {n^2} \sum_{i,j=1}^n |X_n(i,j)|^p < \infty
\]
is finite for all~$p \ge 1$.

The set $\DS_n$ of doubly stochastic matrices of order~$n$ is convex.
The extreme points are the $n!$ permutation matrices,
and the central point $J$ is the equally-weighted average of the extremes.
Given two matrices $A, B \in \DS_n$, the row and column totals of the
deviation matrix $A-B$ are zero.
It is helpful to define the $L^2$-norm and associated metric in $\DS_n$
by the standard Euclidean norm in the tangent space
\[
d^2(A,B) = \sum_{i,j} (A_{ij} - B_{ij})^2,\qquad
\normsq{A} = d^2(A, J) = \sum_{i,j} A_{ij}^2 - 1.
\]
Thus, the central point has norm zero, and each extreme point has
norm $\sqrt{n-1}$.

Each doubly stochastic matrix has one unit eigenvalue,
and all others are less than or equal to one in absolute value.
Consequently, it is helpful to define the operator norm in $\DS_n$ as the Hilbert-Schmidt norm
of the deviation $A-J$, i.e.~the largest absolute eigenvalue $|\lambdamax(A-J)|$.
The spectral gap of $A$ is the difference $\gap(A) = 1 - |\lambdamax(A-J)|$.
%and the boundary of $\DS_n$ consists of those matrices whose spectral gap is zero.
%Thus, the boundary points closest to $J$ are at unit distance.
%In other words, the largest $L^2$-ball contained in $\DS_n$ has unit radius,
%and the smallest $L^2$-ball containing $\DS_n$ has radius $\sqrt{n-1}$.
%The situation is much the same as the cube of side~2 containing
%the unit ball centered at the origin in~$\Real^n$.

A sequence of doubly stochastic matrices $\{A_n\}$ is said to be of
\emph{moderate deviation} if the re-scaled sequence $X_n = n (A_n - J)$ is weakly bounded,
and the spectral gap
$\gap(A_n) \ge C > 0$ is bounded below by a positive constant.
If $A, B$ are moderate-deviation sequences,
the sequence $AB$ of matrix products is also a moderate-deviation sequence,
which implies that $A'A$ is also of moderate deviation.
Convex combinations are also of moderate deviation.
Finally, for fixed~$\nu$, the random sequence $A_n \sim \dsd_n(\nu)$ is of moderate
deviation with probability one.

The main purpose of this note is to show that, for a moderate-deviation sequence
\[
\log \per(n A) = \log(n!) - \vhalf\log \det(I + J - A'A) + O(n^{-1})
\]
with error decreasing as $n\to\infty$.
One consequence is that $\per(nA)/n!$ is bounded as $n\to\infty$.
Even though it does not affect the order of magnitude of the error,
the modified approximation
\begin{equation}\label{det_approx}
\log \per(n A) \simeq \log(n!) - \vhalf\log \det(I + t^2J - t^2 A'A)
\end{equation}
with $t^2=n/(n-1)$ is a worthwhile improvement for numerical work.

\begin{table}
\vbox{
\tabskip 7pt plus10pt
\halign to 12cm{&\hfill$#$ \cr
\noalign{\noindent
Table 1: Numerical illustration of the determinantal approximation.}
\noalign{\smallskip\hrule\medskip}
\omit&\multispan3 \hss $\rho=1$\hss && \multispan3\hss $\rho=2$\hss\cr
n & \log\per(K) & \hbox{Approx} & \hbox{Error}^*&
 & \log\per(K) & \hbox{Approx} & \hbox{Error}^* \cr
\noalign{\smallskip\hrule\medskip}
%5  &3.0397  &3.0250 &0.0147 &&1.8153 &1.7365 &0.0788\cr
%6  &4.5039  &4.4939 &0.0100 &&2.9584 &2.9043 &0.0542\cr
%7  &6.1287  &6.1213 &0.0074 &&4.2820 &4.2425 &0.0395\cr
%8  &7.8905  &7.8846 &0.0059 &&5.7541 &5.7237 &0.0305\cr
%9  &9.7720  &9.7672 &0.0048 &&7.3529 &7.3284 &0.0246\cr
%10 &11.7601 &11.7560 &0.0041 &&9.0627 &9.0422 &0.0205\cr
  8  &7.8905  &7.8895 &1.0178  &&5.7541  &5.7454 &8.7494\cr
 10 &11.7601 &11.7595 &0.6188  &&9.0627  &9.0577 &5.0177\cr
 12 &16.0163 &16.0159 &0.4319 &&12.7701 &12.7667 &3.3501\cr
 14 &20.5955 &20.5951 &0.3293 &&16.8064 &16.8039 &2.4626\cr
 16 &25.4522 &25.4520 &0.2636 &&21.1237 &21.1218 &1.9283\cr
 18 &30.5526 &30.5524 &0.2199 &&25.6868 &25.6852 &1.5778\cr
 20 &35.8700 &35.8698 &0.1880 &&30.4682 &30.4668 &1.3310\cr
 22 &41.3831 &41.3829 &0.1636 &&35.4462 &35.4451 &1.1488\cr
 24 &47.0743 &47.0742 &0.1534 &&40.6032 &40.6022 &1.0134\cr
\noalign{\smallskip\hrule\medskip}
\noalign{${}^*\hbox{Error}\times 10^3$}
}
}
\end{table}

Table~1 illustrates the permanental approximation
applied to a class of symmetric positive definite matrices
$K_{ij} = \exp(-\rho|x_i - x_j|)$ for points $x_1,\ldots, x_n$
equally spaced on the interval $[0, 1]$, so the components
range from $\exp(-\rho)$ to one.
The first step uses the Sinkhorn algorithm to obtain a
doubly stochastic matrix $A$ from~$K$,
and the second step uses the determinantal approximation (\ref{det_approx}).
All values shown are on the log scale.
Because of the difficulty of computing the exact permanent,
the range of $n$-values shown is rather limited.
Nevertheless, it appears that the error decreases 
roughly as $1/(300 n)$ for $\rho=1$, and $1/(40 n)$ for $\rho=2$,
in agreement with~(\ref{det_approx}).
A sequence of matrices of this type with bounded~$\rho$ is of moderate deviation;
a sequence of matrices with $\rho \propto n$ is not.
It is unclear whether a sequence with $\rho=\log(n)$
is moderate or not, but the approximation error does appear to decrease 
roughly at the rate $1/n$ or $\log(n)/n$.

For a very special class of Kronecker-product matrices,
Table~2 shows the relative error in the approximation
for substantially larger values of~$n$.
These matrices have two blocks of size $n/2$ each,
with $B_{ij}=1$ within blocks and $B_{ij}=\rho$ between blocks.
The smaller the between-blocks value, the smaller the spectral gap,
and the greater the approximation error.
For the associated doubly stochastic matrix~$A=2B/(n(1+\rho))$, Table~2 shows
the exact value $\log(\per(nA)/n!)$ together with the error in the
determinantal approximation $-\log\det(I+J-A'A)/2$.
Once again, the empirical evidence suggests that the error
for moderate-deviation matrices with constant~$\rho$ is $O(n^{-1})$, and for
large-deviation matrices, $O(1/(n\rho))$.

\begin{table}
\vbox{%
\offinterlineskip
\tabskip 10pt plus20pt
\halign to \hsize{&\hfill$#$\strut\cr
\noalign{\noindent
Table 2: Determinantal approximation error for large~$n$.\strut}
\noalign{\hrule\smallskip}
\omit&\multispan2 \hss $\rho=0.1$\hss &
	& \multispan2\hss $\rho=0.05$\hss&
	& \multispan2\hss $\rho=5/n$\strut\hss\cr
%\omit &\multispan2\hrulefill&\omit&\multispan2\hrulefill&\omit&\multispan2\hrulefill\cr
n & \hbox{Exact} &  \hbox{Err}\times n& 
  & \hbox{Exact} &  \hbox{Err}\times n&
  & \hbox{Exact} &  \hbox{Err}\cr 
\noalign{\hrule\smallskip}
20 & 0.6755 & 1.6543 &&1.0187 & 2.8616 &&0.3177 & 0.0238 \cr
40 & 0.6274 & 1.3823 &&0.9616 & 3.4398 &&0.5366 & 0.0258 \cr
60 & 0.6140 & 1.2718 &&0.9274 & 3.1088 &&0.6900 & 0.0266 \cr
80 & 0.6082 & 1.2299 &&0.9114 & 2.8589 &&0.8076 & 0.0270 \cr
100 & 0.6049 & 1.2082 &&0.9029 & 2.7229 &&0.9029 & 0.0272 \cr
%120 & 0.6028 & 1.1948 &&0.8977 & 2.6461 &&0.9829 & 0.0274 \cr
%140 & 0.6013 & 1.1856 &&0.8942 & 2.5982 &&1.0519 & 0.0275 \cr
%160 & 0.6002 & 1.1789 &&0.8917 & 2.5653 &&1.1126 & 0.0276 \cr
%180 & 0.5993 & 1.1738 &&0.8898 & 2.5413 &&1.1667 & 0.0277 \cr
200 & 0.5987 & 1.1698 &&0.8882 & 2.5229 &&1.2155 & 0.0278 \cr
%220 & 0.5981 & 1.1666 &&0.8870 & 2.5083 &&1.2600 & 0.0278 \cr
%240 & 0.5977 & 1.1640 &&0.8860 & 2.4965 &&1.3008 & 0.0279 \cr
%260 & 0.5973 & 1.1617 &&0.8852 & 2.4867 &&1.3386 & 0.0279 \cr
%280 & 0.5970 & 1.1598 &&0.8845 & 2.4784 &&1.3737 & 0.0279 \cr
300 & 0.5967 & 1.1582 &&0.8839 & 2.4714 &&1.4065 & 0.0280 \cr
%320 & 0.5964 & 1.1568 &&0.8833 & 2.4653 &&1.4373 & 0.0280 \cr
%340 & 0.5962 & 1.1556 &&0.8829 & 2.4599 &&1.4663 & 0.0280 \cr
%360 & 0.5960 & 1.1545 &&0.8825 & 2.4552 &&1.4937 & 0.0280 \cr
%380 & 0.5958 & 1.1535 &&0.8821 & 2.4511 &&1.5197 & 0.0280 \cr
400 & 0.5957 & 1.1526 &&0.8818 & 2.4473 &&1.5444 & 0.0280 \cr
\noalign{\hrule}
}
}
\end{table}

\section{Justification for the approximation}
\subsection{Permanental expansion}
Let $A$ be a matrix of order~$n$ such that
$n A_i^r = 1 + \epsilon_i^r$, where the row and column totals of $\epsilon$ are zero.
If the components of $\epsilon$ are real numbers greater than $-1$,
$A$~is doubly stochastic and $\epsilon/n = A-J$ is the deviation matrix.
In the calculations that follow, $\epsilon$ is complex-valued and weakly bounded,
and $\epsilon/n$~has spectral norm strictly less than one.
In other words, $A$~need not be doubly stochastic, or even real,
but the associated sequence is assumed to be of moderate deviation.

The permanent expansion of $n A$ as a polynomial of degree~$n$ has the form
\begin{eqnarray*}
\per(1+\epsilon) &=& \sum_\sigma
	(1 + \epsilon_1^{\sigma(1)})
	(1 + \epsilon_2^{\sigma(2)}) \cdots
	(1 + \epsilon_n^{\sigma(n)}) \cr
	&=& n! \Bigl( 1 + n\ave(\epsilon_i^r)
	+ \frac{\dscf n 2}{2!} \avesharp (\epsilon_i^r \epsilon_j^s)
	+ \frac{\dscf n 3}{ 3!} \avesharp (\epsilon_i^r \epsilon_j^s \epsilon_k^t) +
\cdots \Bigr)\cr
	&=& n! \sum_{k=0}^n \frac{\dscf n k}{ k!} \avesharp(\epsilon^{\otimes k}),
\end{eqnarray*}
where $\dscf nk = n(n-1)\cdots(n-k+1)$ is the descending factorial,
and $\avesharp(\epsilon^{\otimes k})$ is the average of $(\dscf n k)^2$ products
$\epsilon_{i_1}^{r_1} \cdots \epsilon_{i_k}^{r_k}$
taken from distinct rows and distinct columns.
More explicitly, the term of degree~$k$ is
$$
%T_k(\epsilon) =
\frac{\dscf n k}{ k!} \avesharp(\epsilon^{\otimes k}) =
	\frac1 {k!\,\dscf nk} \sharpsum \epsilon_{i_1}^{r_1} \cdots \epsilon_{i_k}^{r_k}
$$
with summation over distinct $k$-tuples $(i_1,\ldots, i_k)$ and
$(r_1,\ldots, r_k)$.

In the variational calculations that follow, each term such as the
restricted sum $\sharpsum(\epsilon^{\otimes k})$ or the
unrestricted sum $\sum(\epsilon^{\otimes k})$
is assigned a nominal order of magnitude as a function of~$n$.
For accounting purposes, the components of $\epsilon$ are of order $O(1)$,
and each product $\epsilon_i^r \epsilon_j^s \epsilon_k^t$
is also regarded as being of order one, i.e.~bounded,
at least in a probabilistic sense, as $n\to\infty$.
A sum such as $\sharpsum(\epsilon^{\otimes k})$, for fixed~$k$,
is of order $O((\dscf n k)^2) = O(n^{2k})$,
so the $L^2$-norm $\sum(\epsilon_r^i)^2$ is of order $O(n^2)$.
These sums could be of smaller order under suitable circumstances.
For example, the difference $(\sum - \sharpsum) \epsilon^{\otimes k}$
is a sum of $n^{2k} - (\dscf n k)^2$ terms, and is therefore of order $O(n^{2k-2})$.
Thus, if the sum of the components of $\epsilon$ is zero,
$\sum\epsilon^{\otimes k} = 0$ implies that the restricted sum
$\sharpsum \epsilon^{\otimes k}$ is of order $O(n^{2k-2})$.
Likewise, if the components of $\epsilon$ were independent
random variables of zero mean, $\sum\epsilon_{ir} = O(n)$,
and a similar conclusion holds for the restricted sum.
The variational calculations that follow are not based on
assumptions of statistical independence of components,
but on the arithmetic implications of zero-sum constraints on rows and columns.
For example, $A$~may be symmetric.

The first goal is to show
that the zero-sum restriction on the rows and columns of $\epsilon$
implies that $\sharpsum(\epsilon^{\otimes k})$ is of order
$O(n^k)$ rather than $O(n^{2k})$.
In other words, $\avesharp(\epsilon^{\otimes k})$ is of order $O(n^{-k})$.
In fact, the even-degree terms in the permanent expansion are $O(1)$,
while the odd terms are $O(n^{-1})$.
These conclusions do not hold for all doubly stochastic matrices.
For example, if $A = (\rho I_n + J_n)/(1+\rho)$ for some fixed~$\rho > 0$,
the $L^2$-norm $\norm{A} = \sqrt{n-1}\rho/(1+\rho)$ and other scalars of a similar type
are not bounded. 

\subsection{Zero-sum restriction}
In the term of degree three, the restricted sum is
\begin{eqnarray*}
%n!\, \dscf n k T_3(\epsilon) =
\sharpsum \epsilon_i^r \epsilon_j^s \epsilon_k^t 
	&=& \sharpsum\epsilon_i^r \epsilon_j^s
		(\epsilon_i^r + \epsilon_i^s + \epsilon_j^r + \epsilon_j^s) \cr 
	&=& \sharpsum 2\epsilon_i^r \epsilon_j^s \epsilon_i^r +
		2\epsilon_i^r \epsilon_j^s \epsilon_i^s \cr
	&=& \sharpsum 2(\epsilon_i^r )^3 + 2(\epsilon_i^r)^3 \cr
	&=& 4\sum (\epsilon_i^r)^3. 
\end{eqnarray*}
In the first line, the restricted sum over $k\neq i,j$ is $-\epsilon_i-\epsilon_j$,
while the restricted sum over $t$ of $-\epsilon_i^t$ is
$\epsilon_i^r + \epsilon_i^s$.
Proceeding in this way by restricted summation over each non-repeated index,
we arrive at the following expressions for the restricted tensorial
sums of degree two, three and four:
\begin{eqnarray*}
\sharpsum \epsilon^{\otimes 2} &=&
\sharpsum \epsilon_{i_1, i_2}^{r_1, r_2}  = \sum\epsilon_{i_1,i_1}^{r_1,r_1}
		= \tr(\epsilon'\epsilon) \cr
%\sharpsum \epsilon_{i_1, i_2, i_3}^{r_1, r_2, r_3}  &=
\sharpsum \epsilon^{\otimes 3} &=&
	4\sharpsum \epsilon_{i_1, i_1, i_1}^{r_1, r_1, r_1},  \cr
%\sharpsum \epsilon_{i_1, i_2, i_3, i_4}^{r_1, r_2, r_3, r_4}  &=
\sharpsum \epsilon^{\otimes 4} &=&
	 9\sharpsum \epsilon_{i_1, i_1, i_1, i_1}^{r_1, r_1, r_1, r_1} 
	-9\sharpsum \epsilon_{i_1, i_1, i_1, i_1}^{r_1, r_1, r_2, r_2} 
	-9\sharpsum \epsilon_{i_1, i_1, i_2, i_2}^{r_1, r_1, r_1, r_1} \cr 
	&&{}+ 3\sharpsum \epsilon_{i_1, i_1, i_2, i_2}^{r_1, r_1, r_2, r_2} 
	+6 \sharpsum \epsilon_{i_1, i_1, i_2, i_2}^{r_1, r_2, r_1, r_2} \cr 
	&=& 36\sum \epsilon_{i_1, i_1, i_1, i_1}^{r_1, r_1, r_1, r_1} 
	-18\sum \epsilon_{i_1, i_1, i_1, i_1}^{r_1, r_1, r_2, r_2} 
	-18\sum \epsilon_{i_1, i_1, i_2, i_2}^{r_1, r_1, r_1, r_1} \cr 
	&&{}+3\sum \epsilon_{i_1, i_1, i_2, i_2}^{r_1, r_1, r_2, r_2} 
	+6\sum \epsilon_{i_1, i_1, i_2, i_2}^{r_1, r_2, r_1, r_2}. 
\end{eqnarray*}
Without the zero-sum constraint on the rows and columns,
$\sharpsum \epsilon^{\otimes k}$ is a sum over $(\dscf n k)^2$ $\epsilon$-products,
and thus of order $O(n^{2k})$.
In the reduced form, each distinct value occurs in duplicate at least,
so there are at most $k/2$ distinct values for the row indices and $k/2$ for the column indices.
We observe that $\sharpsum \epsilon^{\otimes 2}=\tr(\epsilon'\epsilon)$ 
is $O(n^2)$, while $\sharpsum\epsilon^{\otimes 3}= 4 \sum(\epsilon_i^r)^3$
is $O(n^2)$ rather than $O(n^3)$.
The five terms in $\sharpsum \epsilon^{\otimes 4}$ are
$O(n^2)$, $O(n^3)$, $O(n^3)$, $O(n^4)$ and $O(n^4)$ respectively.
The final pair can be expressed in matrix notation as
$3\tr^2(\epsilon'\epsilon) + 6 \tr(\epsilon'\epsilon\epsilon'\epsilon)$.

For any vector $x=(x_1,\dots,x_k)$ in $\Real^k$,
let $\tau(x)$ be the associated partition of $[k]=\{1,\ldots, k\}$,
i.e.~$\tau(x)(r,s) = 1$ if $x_r = x_s$ and zero otherwise.
The restricted sum $\sharpsum\epsilon^{\otimes k}$
can be expressed in reduced form either as a restricted sum
or an unrestricted sum,
the only difference arising in the coefficients as shown above for $k=4$.
Generally speaking, unrestricted sums are more convenient for computation.
For general~$k$, the restricted sums are as follows:
\begin{eqnarray*}
\sharpsum \epsilon^{\otimes k}
	&=& \sum_{\rho, \sigma} m^\sharp(\rho) m^\sharp(\sigma)
		\mathop{\sum_{i:\tau(i) = \sigma}}_{r: \tau(r) = \rho} \epsilon_i^r \cr
	&=& \sum_{\rho, \sigma} m(\rho) m(\sigma)
		\mathop{\sum_{i:\tau(i) \ge \sigma}}_{r: \tau(r) \ge \rho} \epsilon_i^r,
\end{eqnarray*}
where the outer sum extends over ordered pairs $\rho, \sigma$ of partitions of
the set $[k]$, and the inner sum over the row and column indices,
$i=(i_1,\ldots, i_k)$ and $r=(r_1,\ldots, r_k)$,
which are held constant in each block. 

The coefficients $m, m^\sharp$ are multiplicative functions
of the partition
$$
m^\sharp(\rho) = \prod_{b\in\rho} (-1)^{\#b-1} (\#b-1),\qquad
m(\rho) = \prod_{b\in\rho} (-1)^{\#b-1} (\#b-1)!,
$$
and $m^\sharp(\rho) = m(\rho)=0$ if $\rho$ has a singleton block.
Here, and elsewhere, $\one_k$ denotes the maximal partition of $[k]$ with one block,
$\#\rho$ is the number of blocks, and $\#b$ is the number of elements of block $b\in\rho$.
Although $m$ and $m^\sharp$ are related by M\"obius inversion,
these expressions are not to be confused with the M\"obius
function for the partition lattice:
$M(\rho, \one) = (-1)^{\#\rho - 1}(\#\rho - 1)!$,
which is a function of the number of blocks independently of their sizes.

For simplicity of notation in what follows, we write
$\epsilon^\rho_\sigma$ for the sum, 
$$
\epsilon_\sigma^\rho =
\mathop{\sum_{i:\tau(i) \ge \sigma}}_{r: \tau(r) \ge \rho} \epsilon_i^r,
$$
in which $i$ and $r$ are constant within blocks,
but the values for distinct blocks may be equal.
Putting these expressions together with a scalar $0 \le t \le 1$, we find
\begin{eqnarray} \label{per_expansion}
\frac{\per(1+t\epsilon)}{n!}
	&=& \sum_{k=0}^n \frac {t^k}{\dscf n k k!}
		\sum_{\sigma,\rho\in\Partition_k} m(\rho) m(\sigma) \epsilon_\sigma^\rho
%\cr
%\nonumber
%	&=& \sum_{k=0}^n \frac {t^k} {k!}
%		\Bigl(1 + \frac{k(k-1)} {2n} + O(n^{-2}) \Bigr)
%		\sum_{\sigma,\rho\in\Partition_k} m(\rho) m(\sigma) \epsilon_\sigma^\rho/n^k \cr
%%	&=& H_n(t) \Bigl(1 + \frac{t^2 H''_n(t)}{2n\, H_n(t)} + O(n^{-2}) \Bigr)\cr
%\nonumber
%	&=& \frac 1 2 \Bigl(H_n(t+t/\surd n) + H_n(t - t/\surd n)\Bigr)(1 + O(n^{-2})),
\end{eqnarray}
where $\Partition_k$ is the set of partitions of~$[k]$.

%If we replace $\dscf n k$ with $n^k$, we obtain the approximation
%%and $H_n(t)$ is obtained by replacing $\dscf n k$ with $n^k$.
%\begin{eqnarray*}
%H_n(t) &=& \sum_{k=0}^n \frac {t^k}{k!}
%                \sum_{\sigma,\rho\in\Partition_k} m(\rho) m(\sigma) \epsilon_\sigma^\rho/n^k \cr
%	&=& \sum_{k=0}^n \frac {t^k} {k!}
%        \sum_{\tau\in\Partition_k} \prod_{b\in\tau} 
%                \sum_{\sigma\vee\rho=\one_b} m(\rho) m(\sigma) \epsilon_\sigma^\rho/n^{\#b} \cr
%h_n(t) = \log H_n(t) &=& \sum_{k=1}^\infty \frac {t^k} {k!}
%                \sum_{\sigma\vee\rho=\one_k} m(\rho) m(\sigma) \epsilon_\sigma^\rho/n^{k},
%\end{eqnarray*}
%where $b \in\tau$ is a block, and $\one_b=\{b\}$ is the maximal one-block partition.
%In effect, $h_n(t)$ is the exponential generator for the
%irreducible matrix-product scalars $\kappa_1,\kappa_2,\ldots$, where
%\[
%\kappa_k =
%\sum_{\sigma\vee\rho=\one_k} m(\rho) m(\sigma) \epsilon_\sigma^\rho/n^{k}
%\]
%is a sum over ordered pairs of partitions $(\rho,\sigma)$ in $\Partition_k$
%whose least upper bound is the full set.

The polynomial (\ref{per_expansion}) is
the exponential generator $\sum_{k=0}^n t^k \zeta_k/k!$ for the moment numbers
\[
\zeta_k = \frac1{\dscf n k} 
	\sum_{\rho,\sigma\in\Partition_k} m(\rho)m(\sigma) \epsilon_\sigma^\rho 
%= \frac 1{\dscf n k} \sharpsum_{\rho,\sigma\in\Partition_k} m^\sharp(\rho)m^\sharp(\sigma)
%\epsilon_\sigma^\rho,
\]
taking $\zeta_k = 0$ for $k > n$.
The coefficients in the log series are the cumulant numbers
\begin{eqnarray*}
\zeta'_k &=& \sum_{\tau\in\Partition_k}
	(-1)^{\#\tau - 1}(\#\tau-1)! \prod_{b\in\tau} \zeta_{\#b} \\
	&=& \sum_{\tau\in\Partition_k}
        (-1)^{\#\tau - 1}(\#\tau-1)! \prod_{b\in\tau} \frac1{\dscf n {\#b}}
\sum_{\rho,\sigma\in\Partition_{b}} m(\rho) m(\sigma) \epsilon_\sigma^\rho \\
	&=& \sum_{\rho,\sigma\in\Partition_k} m(\rho) m(\sigma) \epsilon_\sigma^\rho
	\sum_{\tau\ge \rho\vee\sigma} (-1)^{\#\tau - 1} (\#\tau-1)!
	\prod_{b\in\tau} \frac 1{\dscf n{\#b}} \\
        &=& \sum_{\rho,\sigma\in\Partition_k} m(\rho) m(\sigma) \epsilon_\sigma^\rho
	\times\Delta_n(\rho\vee\sigma)
\end{eqnarray*}
for $k=1,\ldots$.
Note that $\rho,\sigma$ in line~2 are the restrictions to the blocks of $\tau$
of the partitions $\rho,\sigma$ in line~3, so $\tau \ge\rho,\sigma$.
Evidently,
\[
\Delta_n(\rho) = \sum_{\sigma\ge\rho} (-1)^{\#\sigma - 1}(\#\sigma - 1)! \prod_{b\in\sigma}
\frac1{\dscf n {\#b}}
\]
is the generalized cumulant associated with
the reciprocal descending factorial series.
For example, 
\begin{eqnarray*}
\dscf n 4\Delta_n(12|34) &=& 1 - \frac{\dscf n 4} {\dscf n 2 \dscf n 2}
	= \frac {4n-6} {\dscf n 2}\cr
\noalign{\smallskip}
\dscf n 5\Delta_n(123|45) &=& 1 - \frac{\dscf n 5} {\dscf n 3 \dscf n 2}
	= \frac {6(n-2)} {\dscf n 2}\cr
\noalign{\smallskip}
\dscf n 6 \Delta_n(12|34|56) &=& 1 - \frac{3\dscf n 6} {\dscf n 4 \dscf n 2}
		+ \frac{2\dscf n 6}{(\dscf n 2)^3}
	= \frac{8(n-3)(7n-10)}{(\dscf n 2)^2}.
\end{eqnarray*}
Thus, the complete formal expansion for the log permanent is
\begin{eqnarray*}
\log(\per(1+t \epsilon)/n!) 
%	&=& \sum_{k=1}^\infty \frac{t^k \zeta'_k} {k!} \cr 
	&=& \sum_{k=1}^\infty \frac{t^k} {k!}
		\sum_{\rho,\sigma\in\Partition_k} m(\rho)m(\sigma) \epsilon_\sigma^\rho
        \, \Delta_n(\rho\vee\sigma) \cr
	&=& \sum_{k=1}^\infty \frac{t^k} {k!}
		\sum_{\tau \in\Partition_k} \Delta_n(\tau)
		\sum_{\rho\vee \sigma=\tau} m(\rho)m(\sigma) \epsilon_\sigma^\rho
\end{eqnarray*}
in which $\Delta_n(\one_k) = 1/\dscf n k$ for $k\le n$, and
zero otherwise.

\subsection{Moderate-deviation asymptotic expansion}
For a fixed pair of partitions $\rho,\sigma \in\Partition_k$
%whose least upper bound $\rho\vee\sigma$  is the one-block partition $\one_k\in\Partition_k$,
the key scalar $\epsilon_\sigma^\rho$ is the sum
over indices $\tau(i)\ge\sigma$ and $\tau(r)\ge\rho$ of $\epsilon_i^r$.
Under the rules for determining the order of magintude in~$n$,
$\epsilon_\sigma^\rho$
is of order $O(n^{\#\sigma + \#\rho})$, and
\[
\Delta_n(\rho\vee\sigma)\, \epsilon_\sigma^\rho = O(n^{\#\sigma+\#\rho - k - \#(\rho\vee\sigma) + 1}).
\]
Since neither partition has singleton blocks, $2\#\sigma \le k$
and $2\#\rho\le k$.
Thus the order is $O(1)$ only if $k$ is even,
each partition is binary with all blocks of size two,
and the least upper bound is the full set~$\one_k$.
Otherwise, if the least upper bound is less than $\one_k$,
or if any block of either partition is of size three or more,
the term is $O(n^{-1})$ or smaller.

%By collecting irreducible scalars of equal order in~$n$, we may write
%$h_n(t)$ as an asymptotic expansion $h_n^{(0)}(t) + h_n^{(1)}(t) + \cdots$ where
%\begin{eqnarray*}
%h_n^{(r)}(t) &=& 
%	\sum_{k=1}^\infty \frac {t^k} {k!}
%                \mathop{\sum_{\sigma\vee\rho=\one_k}}_{\#\sigma+\#\rho=k-r}
%		 m(\rho) m(\sigma) \epsilon_\sigma^\rho/n^{k}
%\end{eqnarray*}
%is the sum of irreducible scalars of order $O(n^{-r})$.
%The outer sum can be restricted to terms of degree $k\ge r+2$.
%
The leading term of maximal order in the expansion of the log permanent is
\begin{eqnarray*}
h_n^{(0)}(t) 
	&=& \sum_{k=1}^{n/2} \frac{t^{2k}} {(2k)!\, \dscf n{2k}}
	\mathop{\sum_{\rho,\sigma \sim 2^k}}_{\rho\vee\sigma=\one_{2k}} \epsilon_\sigma^\rho \cr 
	&=& \sum_{k=1}^{n/2} \frac {t^{2k}} {2k} \tr((\epsilon'\epsilon)^k)/\dscf n{2k} \cr
	&=& \sum_{k=1}^\infty \frac {t^{2k}} {2k} \tr((\epsilon'\epsilon/\dscf n 2)^k) + O(n^{-1}) \cr
\noalign{\smallskip}
	&=& -\vhalf\log\det(I_n - t^2\epsilon'\epsilon/\dscf n2) + O(n^{-1}) \cr
\noalign{\smallskip}
	h_n^{(0)}(1) &=& -\vhalf\log\det(I_n + J - A'A) + O(n^{-1}).
\end{eqnarray*}
The symbol $\rho,\sigma\sim 2^k$ denotes two partitions of $[2k]$ having $k$
blocks of size two, and
there are $(2k-1)!$ pairs whose least upper bound is the full set $[2k]$.
The series is convergent for $|t| < 1/|\lambdamax(A-J)|$,
so $h_n^{(0)}(1)$ is finite by the spectral gap assumption.

%In the passage from $\dscf n {2k}$ to $(\dscf n 2)^k$, certain terms of order $O(n^{-1})$
%have been dropped.
%Since $1/\dscf n {2k} - 1/n^{2k} \simeq k(2k-1)/n^{2k+1}$
%the error is roughly $\sum (k-1/2) \lambda_i^k / n$, where $\{\lambda_i\}$ are the eigenvalues
%of $t^2\epsilon'\epsilon/n^2$.
%The correction is slightly different if $\epsilon'\epsilon/n^2$ were replaced
%with $\epsilon'\epsilon/\dscf n 2$ in lines 3--4, as suggested in section~2.

In the expansion of the first-order correction $h^{(1)}(t)$,
the terms may be grouped by degree in~$t$.
The following are the terms of degree six or less,
expressed so far as possible using matrix operations in which $\epsilon_r$
and $(\epsilon'\epsilon)_r$ are
the component-wise $r$th powers of $\epsilon$ and $(\epsilon'\epsilon)$
respectively.
\begin{eqnarray*}
\hbox{degree 3:\quad}&& \frac 2 {3n^3} \sum \epsilon_3 = \frac 2 {3n^3}\sum (\epsilon_i^r)^3 \cr
\noalign{\vskip 4pt}
\hbox{degree 4:\quad}&& -\frac{3}{4n^4} \sum (\epsilon_2'\epsilon_2 + \epsilon_2\epsilon_2') \cr
\noalign{\vskip 4pt}
\hbox{degree 5:\quad}&& \frac 2{n^5}\tr(\epsilon_2 \epsilon'\epsilon\epsilon')
	+ \frac 1 {n^5} \sum \epsilon_2\epsilon' \epsilon_2 \cr
\noalign{\vskip 4pt}
\hbox{degree 6:\quad}&& \frac 1 {3n^6} \sum ((\epsilon'\epsilon)_3 + (\epsilon\epsilon')_3) 
	+ \frac 1 {2n^6} \sum(\epsilon'_2 \epsilon\epsilon'\epsilon_2
		+ \epsilon_2 \epsilon'\epsilon\epsilon'_2) \cr
\noalign{\vskip 4pt}
	&&\quad -\frac 3 {2n^6} \sum (\epsilon_2(\epsilon'\epsilon)_2 +
	\epsilon'_2 (\epsilon\epsilon')_2 )
\end{eqnarray*}

\section{Random matrices}
In order to test the adequacy of the determinantal approximation,
5000~random matrices $A$ of order 10--25 were generated,
and the exact values $y(A) =\log(\per(n A)/n!)$ were computed using the Ryser algorithm.
The matrices of even order were generated from the distribution $A_i \sim \dsd_n(\nu_i)$.
%It was observed that $2y(A) \simeq 1/\hat\nu(A)$, with a correlation of around~0.5
%for fixed~$\nu$.
In order to ensure an adequate range of permanental values, the parameters $\nu_i$ were
generated randomly and independently from the Gamma distribution
with mean one on four degrees of freedom, i.e.~$\chi_8^2/8$.
On the matrices of odd order was superimposed a multiplicative
random block pattern $(\eta + \eta' B(r,s))$
in which $\eta,\eta'$~are uniform $(0,1)$ random scalars,
and $B$ is a random partition of $[n]$, the Chinese restaurant process with parameter~1.
This block pattern was superimposed prior to Sinkhorn projection in
order to vary the eigenvalue pattern, which depends on the coefficients~$\eta$
and on the block sizes.
Matrices of this type do not satisfy the moderate-deviation condition.

\begin{figure}
\vbox{\hsize 13.5cm
\includegraphics[width=13.5cm, height=17cm]{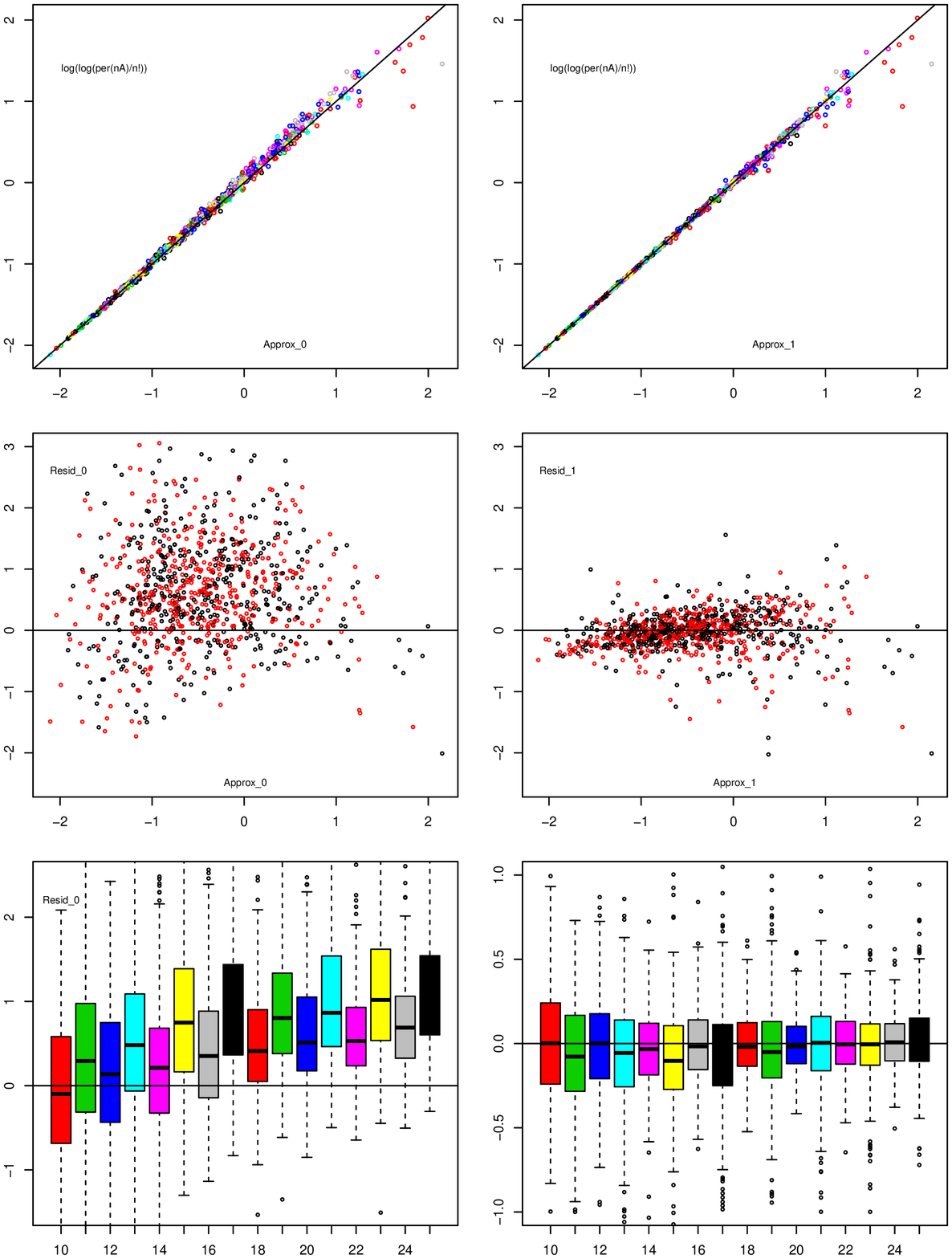}
\vskip 1pc
\vbox{%
\noindent
Figure 1. Top row: exact permanent of 500 random matrices plotted
against two approximations, both on the log-log scale.
Middle row: standardized log-log residuals plotted against the two approximations.
Third row: residuals plotted against~$n$, at twice the scale in the right panel.
}}
\end{figure}
 
On the log scale, the range of observed permanental values was $0.11 < y(A) < 7.56$,
or $0.042 < y/\log(n) < 2.57$, the largest value occurring for a matrix of order~19.
The target range $y(A) \le 2$ was exceeded by 2.6\% of the matrices generated,
and the moderate-deviation threshold $y > \log(n)$ was exceeded by
43 matrices comprising roughly 0.9\% of the simulations.
Most of these exceedances occurred for matrices of odd order having a pronounced
block structure.
The range of $L^2$-values was $0.22 \le \normsq{A}\le 7.32$,
with only~2\% in excess of~3.0.
Although there exist matrices $A\in\DS_n$ such that the ratio $y(A) / \normsq{A}$ is
arbitrarily large or arbitrarily small, the range of simulated values was
only $(0.53, 1.15)$.
This is a reminder that the behaviour of the permanent in a bounded $L^2$-ball
is not a good indicator of its behaviour near the corners.

Two approximations were computed,
the first-order determinantal approximation $x_0(A)$ as in~(\ref{det_approx}),
and a second-order modification $x_1(A)$ using the additional terms in section~4.3.
The scatterplots of $\log\log y(A)$ against $\log\log x_0(A)$ and $\log\log x_1(A)$
are shown in the top two panels of Fig.~1.
To reduce clutter in the top two rows, only a 10\% sample is shown.

Since the relative error of the approximation $x(A)$ increases with its magnitude,
the standardized residuals are most naturally defined on the log-log scale
\[
\hbox{Resid}(y, x) = \frac{n\,\log(y(A)/x(A))} {x(A)}.
\]
In the middle panels of Fig.~1, the log-log residuals are plotted against~$x$,
using the same scale for both plots.
In the lower panels, the log-log residuals are plotted against~$n$.
The first plot makes it clear that the determinantal approximation $x_0(A)$
tends to be an under-estimate:
the rate of occurrence of the inequality $y(A) > x_0(A)$ increases from
44\% for $n=10$ to over 90\% for $n\ge 20$.
In the lower left panel, the alternating pattern of residuals for odd and even~$n$
is a consequence of the block pattern embedded in the matrices of odd order.
For both types of matrices, the plots suggest that the residual distribution
is asymptotically constant as $n\to\infty$,
and that the correction term is appreciable for moderate~$n$.
The root-mean-squared residual is approximately $1.2/n$ for the first approximation
and $0.3/n$ for the second, but the distributions are more like Cauchy than Gaussian.
As it happens, the inequalities
%resid0 <- nn*log(y/x0)/x0       # 99% point around -1.4, 2.5
%resid1 <- nn*log(y/x1)/x1       # 99% point around \pm 0.6
\[
\log y(A) < \log x_0(A) - 1.4\,x_0(A)/n,\qquad
\log y(A) > \log x_0(A) + 3.0\,x_0(A)/n
\]
occur in the sample with rates 1.0\% each.
% respectively for $n\ge 20$.
The inequality
\[
|\log y(A) - \log x_1(A)| > 0.5 x_1(A)/n
\]
occurs at a fairly constant rate of around 1\% for non-block-structured matrices.
For block-structured matrices, the rate is about five times as high,
but slowly decreasing in~$n$ over the range observed.
Of course, these rates must depend on the distribution by which the
matrices are generated,
but it is anticipated that the rate will have a positive limit for matrices
in the moderate-deviation region $x_0(A) < \log(n)$.
In other words, we should expect $\log(\per(nA)/n!)$ to lie in the
interval $x_1(A) \exp(\pm x_1(A)\log(n)/n)$ with probability
tending to one for large~$n$ if $x_1(A)$ is bounded.

%If $A$ is the Sinkhorn projection of a matrix $Y$ whose components are
%independent and uniformly distributed, the approximation $x_1(A)$ is accurate to about
%three significant decimal digits for $10\le n \le 25$.
%If the components are independent Bernoulli variables with parameter $1/2$,
%the approximation $x_1(A)$ is accurate to about ???.

\parskip 0.5pc \parindent 0pt \hangindent 10pt
\section{References}
Chatterjee, S., Diaconis, D. and Sly, A. (2010)
Properties of uniform doubly stochastic matrices.
ArXiv:1010.6136v1 [Math.PR]

Deming, W.E. and Stephan, F.F. (1940)
On a least-squares adjustment of a sampled frequency table when the expected
marginal totals are known.
{\it Ann.\ Math.\ Statist.} 11, 427-444.

Jerrum, M. and Sinclair, A. (1989)
Approximating the permanent.
{\it SIAM J.~Comput,}, 18, 1149--1178.

Jerrum, M., Sinclair, A. and Vigoda, E. (2004)
A polynomial-time approximation algorithm for the permanent of a
matrix with non-negative entries.
{\it Journal of the ACM} 51, 671-697.

Kou, S. and McCullagh, P. (2009).
Approximating the $\alpha$-permanent.
{\it Biometrika\/} 96, 635--644.

Linial, N., Samordnitsky, A. and Wigderson, A. (1998)
A deterministic strongly polynomial algorithm for matrix scaling
and approximate permanents.
{\it Proc.\ 30th ACM Symp.\ on Theory of Computing.}
ACM, New York.

McCullagh, P. and Nelder, J.A. (1989)
{\it Generalized Linear Models.}
London, Chapman \& Hall.

Sinkhorn, R. (1964)
A relationship between arbitrary positive matrices and 
doubly stochastic matrices.
{\it Ann.\ Math.\ Statist.} 35, 876--879.

Valiant, L. (1979)
The complexity of computing the permanent.
{\it Theoretical Computer Science} 8, 189--201.

\end{document}